\pdfoutput=1
\RequirePackage{ifpdf}
\ifpdf 
\documentclass[pdftex]{sigma}
\else
\documentclass{sigma}
\fi

\DeclareMathOperator{\const}{const}

\begin{document}

\newcommand{\dd}[2] {{{\partial #1} \over {\partial #2}}}

\newcommand{\arXivNumber}{1407.7163}

\allowdisplaybreaks

\renewcommand{\PaperNumber}{101}

\FirstPageHeading

\ShortArticleName{Who's Afraid of the Hill Boundary?}

\ArticleName{Who's Afraid of the Hill Boundary?}

\Author{Richard MONTGOMERY}

\AuthorNameForHeading{R.~Montgomery}

\Address{Math Dept.
UC Santa Cruz, Santa Cruz, CA 95064, USA} \Email{\href{rmont@ucsc.edu}{rmont@ucsc.edu}}
\URLaddress{\url{http://people.ucsc.edu/~rmont/}}

\ArticleDates{Received August 25, 2014, in f\/inal form October 28, 2014; Published online November 02, 2014}

\Abstract{The Jacobi--Maupertuis metric allows one to reformulate Newton's equations as geodesic equations for
a~Riemannian metric which degenerates at the {\it Hill boundary}.
We prove that a~JM geodesic which comes suf\/f\/iciently close to a~regular point of the boundary contains pairs of
conjugate points close to the boundary.
We prove the conjugate locus of any point near enough to the boundary is a~hypersurface tangent to the boundary.
Our method of proof is to reduce analysis of geodesics near the boundary to that of solutions to Newton's equations in
the simplest model case: a~constant force.
This model case is equivalent to the beginning physics problem of throwing balls upward from a~f\/ixed point at f\/ixed
speeds and describing the resulting arcs,
see Fig.~\ref{envelope}.}


\Classification{37J50; 58E10; 70H99; 37J45; 53B50}

\section{Results and motivation}
One constructs the Jacobi metric $ds^2_{\rm JM}$ of classical mechanics by f\/ixing the total energy~$E$ of the system and
multiplying the kinetic energy metric $ds^2_{\rm K}$ by the conformal factor
\begin{gather*}
f =2(E-V),
\qquad
ds^2_{\rm JM} = f ds^2_{\rm K},
\end{gather*}
where~$V$ is the potential energy.
It is well-known that the geodesics for this Jacobi--Maupertuis metric (henceforth JM metric for short) are, up to
reparameterization, exactly the solutions to Newton's equations having energy~$E$.
(See Proposition~\ref{JMtoNewton} below for a~careful statement.
See~\cite[Theorem 3.7.7]{AbMar} for another discussion and a~nice proof.) The domain of the Jacobi metric is the domain
in conf\/iguration space where this conformal factor is non-negative and is called the Hill region:
\begin{gather*}
{\mathcal H} = \{q: f (q)\ge 0\}.
\end{gather*}
The Hill region includes the Hill boundary (sometimes called the zero velocity surface) where the conformal factor, and
hence the metric, vanishes:
\begin{gather*}
\partial {\mathcal H}= \{q:f(q) = 0\}.
\end{gather*}
A~``regular point'' $q_0$ of the Hill boundary is one for which $df(q_0) \ne 0$.
Here is our main result.
\begin{theorem}\label{Theorem1}
Any JM geodesic which comes sufficiently close to a~regular point $q_0$ of the Hill boundary contains a~pair of
conjugate points close to $q_0$ which are conjugate along a~short arc close to $q_0$.
In particular, such a~geodesic fails to minimize JM length.
\end{theorem}

This theorem is a~direct consequence of a~structure theorem, Theorem~\ref{Structure} below, regarding the conjugate
locus of near-boundary points, and results from Seifert's seminal paper~\cite{Seifert} which we recall in the next
section.

{\bf Motivations.} Two questions motivated this paper.

1.~Can the calculus of variations, applied to the JM metric reformulation of mechanics, uncover new results regarding the
classical three-body problem? The direct method of the calculus of variations breaks down at the Hill boundary since
curves lying in the boundary have zero JM length.
A~deeper understanding of the behaviour of near-boundary JM geodesics seems necessary to the further development of JM
variational methods in case where the Hill boundary is not empty.
For some results in celestial mechanics based on JM variational methods in instances where the Hill boundary is not
empty see~\cite{Moeckel} and~\cite{Soave_Terr} whin this direction

2.~Does the fact that JM curvatures tend to positive inf\/inity imply there are conjugate points near the boundary? Let~$q$
be a~point near a~regular point of the Hill boundary and let~$y$ denote its Riemannian distance from the boundary.
The sectional curvatures~$K$ of two-plane through~$q$ which contains the normal direction to the boundary tends to
positive inf\/inity like $1/ y^3$ as $y \to 0$.
The classical Bonnet--Meyer's estimate says that if the curvatures~$K$ along a~geodesic through~$q$ are greater than or
equal to a~positive constant $K_0$ then there must be a~point conjugate to~$q$ along the geodesic and lying within
$\pi/\sqrt{K_0} $ from~$q$.
This suggests the existence of conjugate points within $y^{3/2}$ from our point~$q$.
However, the {\it JM distance} of~$q$ to the boundary is also of order $y^{3/2}$ for small~$y$.
The two distances are of the same order.
These naive estimates do not tell us if Bonnet--Meyers ``wins'' to beat out the closeness of the boundary by creating
a~conjugate point before we have ``ref\/lected'' of\/f the boundary and left the region in which the Bonnet--Meyers curvature
estimate holds.
Theorem~\ref{Theorem1} asserts that, indeed, Bonnet--Meyers wins.

3.~The recent work~\cite{Ramis} claims that the harmonic oscillator, when it is reformulated in terms of JM geodesics, has
positive Lyapunov exponents.
This surprise, and trying to better understand it, was the seed that planted this paper.

\section{Mechanics and Seifert's coordinates}

By Newton's equations on a~manifold~$M$ we mean a~system of second-order dif\/ferential equations of the form
\begin{gather}
\nabla_{\dot \gamma} \dot \gamma = - \nabla V (\gamma).
\label{N}
\end{gather}
Here $\nabla$ is the Levi-Civita connection associated with a~f\/ixed Riemannian metric $ds^2_{\rm K}$ on~$M$.
(The subscript `K' is for `kinetic'.)
$V$ is a~chosen smooth function on~$M$ called the ``potential''.
The total energy
\begin{gather*}
H = \frac{1}{2} \langle \dot \gamma, \dot \gamma \rangle_{\gamma} + V(\gamma)
\end{gather*}
is constant along any solution to~\eqref{N}.
The inner product is the one def\/ined by the metric $ds^2_{\rm K}$.

Fix a~value $H=E$ for this energy and form the conformal factor
\begin{gather*}
f = 2(E-V)
\end{gather*}
and the resultant Jacobi--Maupertuis metric
\begin{gather}
ds^2_{\rm JM} = f ds^2_{\rm K}.
\label{JM}
\end{gather}
The following well-known proposition connects solutions to~\eqref{N} with Jacobi geodesics.
\begin{proposition}\label{JMtoNewton}
Solutions to~\eqref{N} with energy~$E$ are, after reparameterization, geodesics for the metric~\eqref{JM}  
which
lie inside the Hill region $ f \ge 0$ and touch the Hill boundary $f=0$ in at most two points.
Conversely, any geodesic for the Jacobi metric lying inside the Hill region and touching the boundary in no more than
two points is a~reparameterization of a~solution to~\eqref{N}.
\end{proposition}

For a~proof see~\cite[Theorem 3.7.7]{AbMar}.

Special care must be taken with geodesics at the Hill boundary.
We have $f = 2(E -V(\gamma(t)) = \| \dot \gamma \|^2)$
along solutions $\gamma(t)$ to~\eqref{N}.
It follows that such a~solution hits the boundary at a~time $t_0$ if and only if $\dot \gamma (t_0) = 0$.
We call such a~solution a~``brake orbit''.
The point $q_0 = \gamma (t_0)$ where the solution hits the boundary is called the ``brake point'' since it has
instantaneously stopped.
Uniqueness of solutions to~\eqref{N} shows that a~brake orbit retraces its own path when we pass the brake instant:
$\gamma(t_0 + t) = \gamma(t_0 -t)$.
When we speak of Jacobi geodesics which hit the Hill boundary we mean exactly these brake orbits, up to
reparameterization.

If a~brake orbit hits the Hill boundary at two distinct points then it is periodic, shuttling back and forth forever
between these two brake points, with Newtonian period twice the Newtonian time it takes to get from one point to the
other.
Conversely, any periodic orbit having one brake point must have another distinct brake point.
Seifert's primary aim in~\cite{Seifert} was to establish the existence of such periodic brake solutions.

Suppose that the brake point $q_0$ is a~regular point of the boundary: $df(q_0) \ne 0$, i.e.~$\nabla V(q_0) \ne 0$.
Then Seifert proved that for small~$\epsilon$ the sub-arc $\gamma([t_0, t_0 + \epsilon])$ of~$\gamma$ is a~minimizing JM
geodesic which realize the JM distance from $\gamma(t_0 + \epsilon)$ to the boundary.
A~Taylor expansion yields $\gamma(t_0 + h) = q_0 - \frac{1}{2} h^2 \nabla V(q_0) + O(h^4)$ showing that this brake
orbit, as a~non-parameterized curve, is smooth and intersects the boundary orthogonally at the brake point.

Seifert solved Newton's equations with initial conditions $(\gamma, \dot \gamma) = (q, 0)$ on the Hill boundary to form
a~system of coordinates $(x, y)$ with $x = (x_1, \ldots, x_{n-1})$, $x_i, y \in \mathbb{R}$ and $n=\dim(M)$ for which
the~$y$-curves $x_i=\const$, $y=t$ are reparameterized brake orbits with brake instant $t=0$.
In these coordinates the Hill boundary is given by $y=0$ and the~$x_i$ coordinatize points on the Hill boundary.
We center the coordinates at a~regular point $q_0$ of the boundary, meaning that~$q_0$ has coordinates $(0,0)$.

We will call such coordinates ``cylinder coordinates'' or ``Seifert coordinates''.
When the coordinate domain has the form $W \times [0, \epsilon]$ we call the resulting sets in the manifold ``cylinder
neighborhoods'' of height~$\epsilon$.
The ``roof'' of the cylinder is the locus $y = \epsilon$.
The ``vertical lines'' are the images of $\{q_* \} \times [0, \epsilon]$ and are brake orbits.
We will say that the direction $-\dd{}{y}$ is `straight down'.
It is the tangent direction f\/ield to the brake orbits headed to the boundary.

{\bf Properties of Seifert coordinates.}
We recall some properties of Seifert's cylinder coordinates.
Within the f\/ive bullet points
\begin{itemize}\itemsep=0pt
\item[(1)] The JM distance of a~point to the Hill boundary is $y^{3/2}$
(\cite{Seifert}, the displayed equation just preceding his equation (46)).
\item[(2)] The JM metric is given by the coordinate expression
\begin{gather}
ds^2_{\rm JM} = ydy^2 + yf(x,y)\Big(\sum dx_i^2 + \sum h_{ij} (x,y) dx_i dx_j\Big),
\label{SeifertModel}
\end{gather}
where $f(0,0) = 1$ and $h_{ij}(x,y) = O(x^2+y^2)$ \cite[equation~(46)]{Seifert}.
\item[(3)]
For any suf\/f\/iciently small cylinder neighborhood~$A$ of $q_0$ and any~$\delta$ smaller than $45$ degrees
there is a~smaller cylinder neighborhood~$B$ of $q_0$ such that every geodesic which enters into~$B$ must exit and
leave~$A$ through the roof of~$A$, making an angle of less than~$\delta$ with the vertical line as it enters and leaves
(see Fig.~\ref{parabola}).

\begin{figure}[t]\centering
\includegraphics[width=5.0cm]{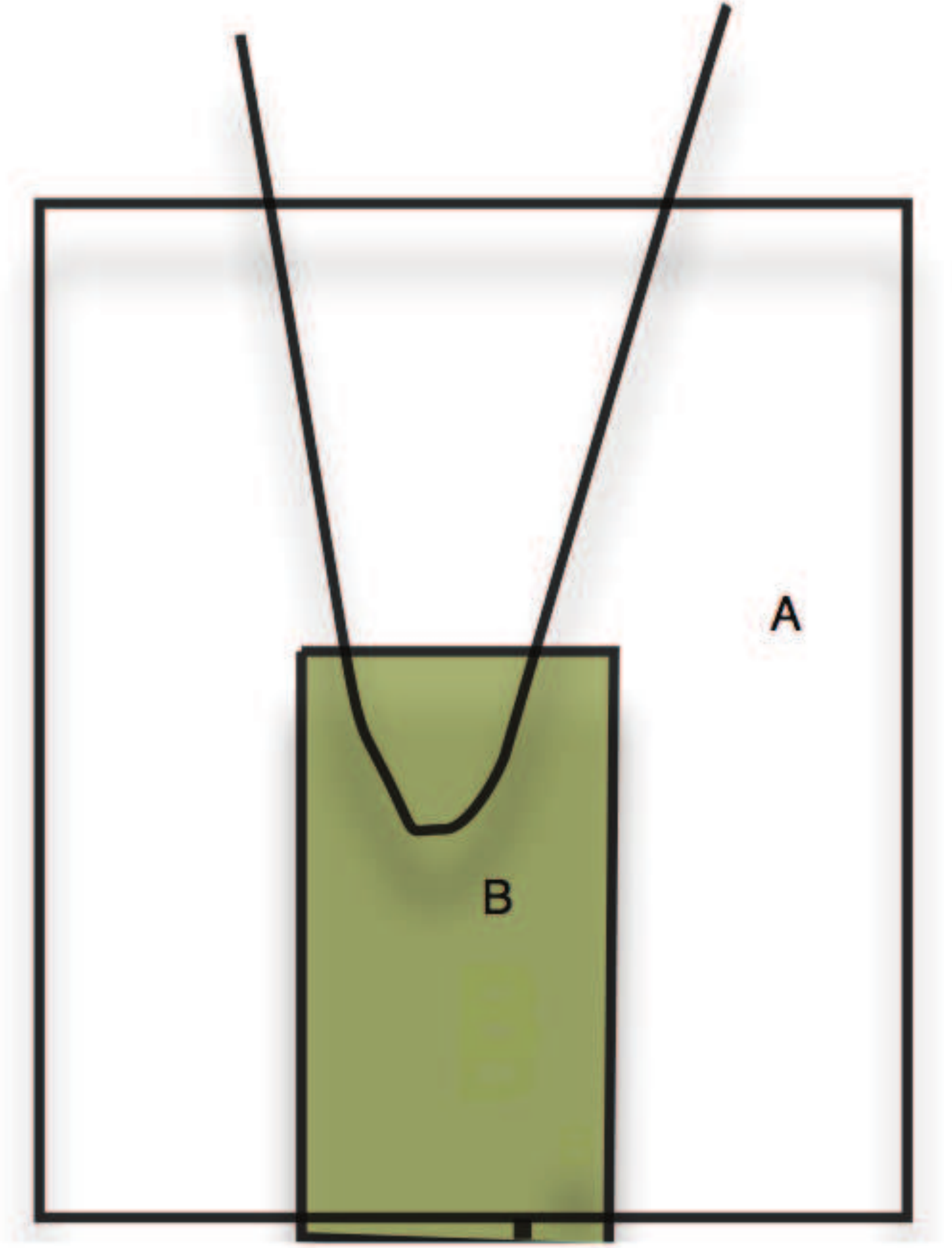}
 \caption{A geodesic which enters into $B$ must enter and leave
through the roof of $A$~at a~steep upward angle.}
\label{parabola}
\end{figure}

\item[(4)]
Along any of the geodesics $\gamma(t) = (x(t), y(t))$ described in the previous item, the height function
$y(t)$ is strictly convex relative to Newtonian time~$t$, with a~unique local minimum~\cite[Fig.~3, also Theorem~\ref{Theorem1}]{Seifert}.

\item[(5)] If a~geodesic enters into a~suf\/f\/iciently small cylinder neighborhood of height~$h$ then it leaves that
neighborhood within a~short Euclidean time of at most $C \sqrt{h}$ where~$C$ is any constant greater than $2 \sqrt{2}/\|\nabla V (q_0)\|$ \cite[Fig.~3 and Theorem~\ref{Theorem1}]{Seifert}.
\end{itemize}

Items (3) and (4) are not proved exactly as stated in Seifert.
We give proofs in Section~\ref{Section5} below.

For~$q$ a~point in a~cylinder set~$A$ let $q_* \in \partial {\mathcal H} \cap A$ denote the brake point along the brake
orbit connecting~$q$ to the boundary.
In terms of cylinder coordinates, if $q = (x,y)$ then $q_* = (x,0)$.
Let $C(q) \subset \bar {\mathcal H} \cap A$ denote the f\/irst conjugate locus to~$q$ for the restriction of the JM metric
to~$A$.
The points of $C(q)$ are the points conjugate to~$q$ along geodesic arcs {\it lying in}~$A$.

\begin{theorem}[structure theorem]\label{Structure}
Let~$A$ be a~cylinder set whose height~$\epsilon$ is sufficiently small.
Then the conjugate locus $C(q) \subset A$ of any point $q \in A$ has the following properties.
$C(q)$ is a~smooth hypersurface which intersects the Hill boundary tangentially at $q_*$ and in no other point.
As a~singularity of the exponential map, $C(q)$ represents the fold singularity.
Every geodesic arc through~$q$ in~$A$ lies entirely on the side of $C(q)$ closest to~$q$.
With the single exception of the brake orbit $[q, q_*]$, if such a~geodesic arc touches $C(q)$ then it touches it
tangentially.
Every geodesic through~$q$ whose initial tangent vector~$v$ is sufficiently close to the ``straight down'' direction
touches $C(q)$ $($see Fig.~{\rm \ref{envelope})}.
\end{theorem}

\begin{remark}
 Compare this theorem with part (C) of Theorem 3.3 of~\cite{Warner} where Warner shows fold-type
conjugate loci occur stably and generically for Riemannian metrics.
\end{remark}

\begin{remark} If we take a~geodesic which starts at~$q$ and touches $C(q)$ and extend it slightly beyond $C(q)$ then
it will fail to minimize.
The extent to which it fails to minimize is measured by the {\it index} of Morse theory.
This index is 1 for all the geodesics of the structure theorem, this being the dimension of the kernel described towards
the end of Appendix~\ref{A2} in the paragraph {\bf Fold}.
\end{remark}

Let us continue with the notation of Theorem~\ref{Structure}.
If~$v$ is a~tangent vector to $q \in A$ then we say that~$v$ ``points downward'' if $dy(v) < 0$.
Consider the cone of downward-pointed velocities $v \in T_q {\mathcal H}$ with the additional property that the geodesic
with initial condition $(q,v)$ touches the conjugate locus $C(q)$ to~$q$ at a~point $c \in C(q)$ below~$q$: $y(c) <
y(q)$.
Call this set of vectors~$v$ the ``downward conjugate cone'' at~$q$ and denote it by ${\rm DC}(q)$.

\begin{theorem}\label{Theorem3}
The downward conjugate cone ${\rm DC}(q)$ is an open cone containing the brake direction.
As $q \to q_*$ along the brake segment $[q, q_*]$, the cone ${\rm DC}(q)$ limits to the open downward pointed cone consisting
of all vectors~$v$ whose angle with the straight down direction is less than $45$~degrees.
\end{theorem}

\section{Throwing balls: the model example}

The idea of our proof is to reduce the study of geodesics near the Hill boundary to that of a~model example for which
the geodesics can be found exactly.
The model is
\begin{gather}
\label{model}
ds^2_{\text{falling}} = y\big(dx^2 + dy^2\big),
\qquad
y\ge 0.
\end{gather}
In the model $f = y$, $ds^2_{\rm K} = dx^2 + dy^2$, $V(x,y) = -\frac{1}{2}y$ and $E =0$.
The corresponding Newton's equations are
\begin{gather}
\ddot x = 0,
\qquad
\ddot y = 1/2
\label{modelEqs}
\end{gather}
with energy
\begin{gather*}
H(x,y,\dot x, \dot y) = \frac{1}{2}\big(\dot x^2 + \dot y^2\big) -\frac{1}{2}y.
\end{gather*}

\begin{figure}[h]\centering \includegraphics[width=65mm]{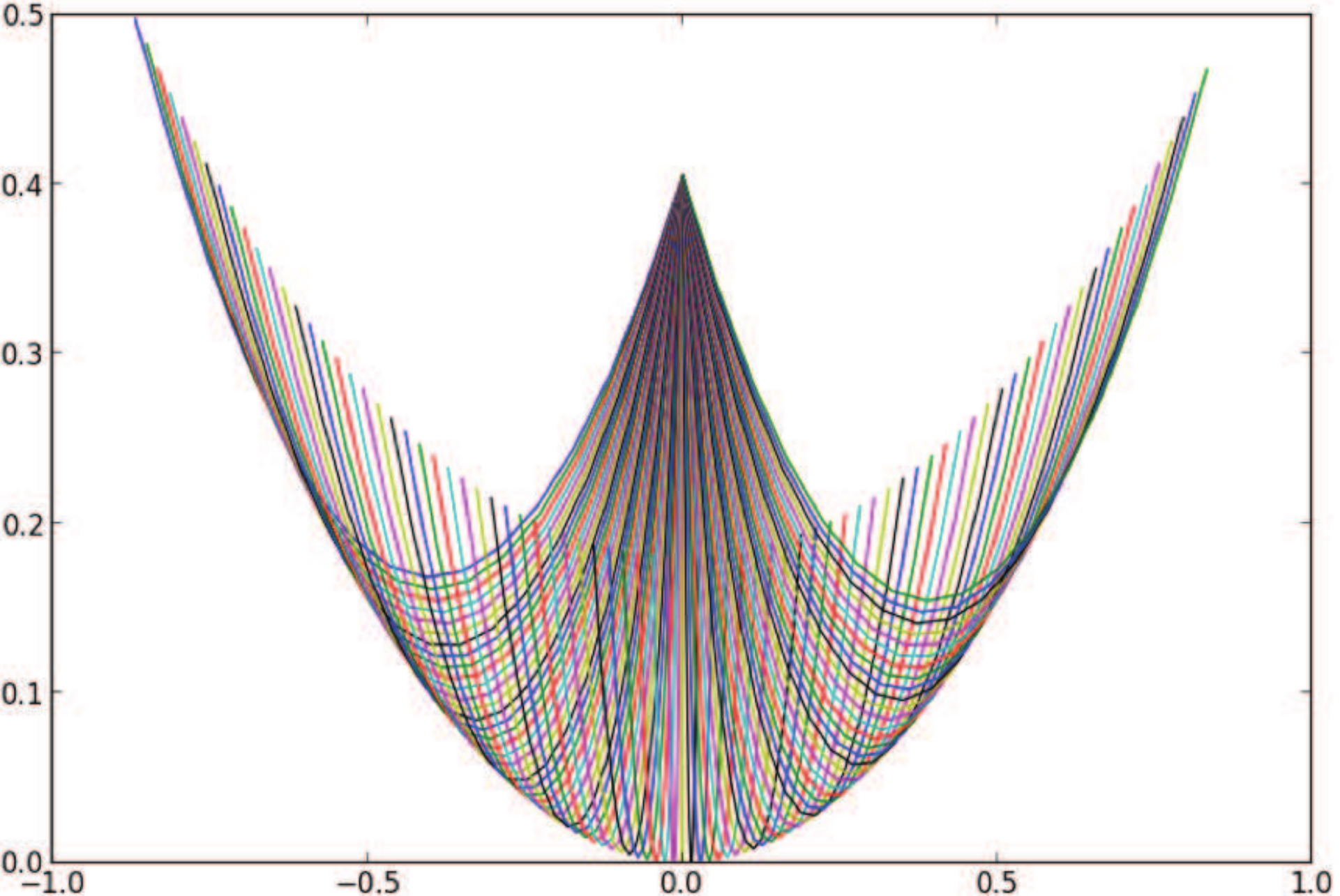}
 \caption{Geodesics leaving a~point and headed toward the
boundary in the model example.
Turn the f\/igure upside down to see the trace of thrown balls, a~sprinkler, or f\/ireworks.}
\label{envelope}
\end{figure}

{\bf Freshman physics.} The af\/f\/ine change of variables $z = h_0 - 2gy$, $x =x$
turns these Newton's equations into the equation $\ddot z = -g$, $\ddot x = 0$
which governs the height~$z$ of a~ball thrown under the inf\/luence of the earth's
constant gravitational f\/ield of strength~$g$, pointed down.
This is the well-studied problem of ballistics from the 1st week or so of most beginning physics courses.
We are throwing balls or shooting cannons from a~f\/ixed point~$q$ with $z(q) < h_0$.
The Hill region $f \ge 0$ is $z \le h_0$.
The speed of our throws at a~f\/ixed point are all equal and are such that the maximum possible height we can reach, the
height reached if we hurl our ball straight up, is the height~$h_0$.
Turn Fig.~\ref{envelope}
upside down to see a~familiar picture of many balls being thrown at the same speed from
a~f\/ixed point to form a~sprinkler pattern, or if you prefer, the arcs of light seen in a~f\/ireworks display.

The general solution to our model Newton's equations~\eqref{modelEqs} is the family of parabolas:
\begin{gather}
x = x_0 + v_1 t,
\qquad
y = y_0 + v_2 t - \frac{1}{4} t^2.
\label{Solutions}
\end{gather}
The parameters $(x_0,y_0, v_1, v_2)$ are the initial conditions at time $t=0$ for our dif\/ferential
equations~\eqref{modelEqs}.
The energy along any one member of this family of solutions is $H(x_0,y_0, v_1, v_2)= \frac{1}{2} (v_1^2 + v_2^2)
-\frac{1}{2}y_0$.
We want this energy to be zero which means that $v_1^2 + v_2^2 = y_0$ so that the allowable velocities $(v_1, v_2)$
through $P_0$ vary over a~circle.
For each velocity in this circle we get a~parabola through $P_0$.
The {\it envelope} of this one-parameter family of parabolas is the conjugate locus.

\begin{lemma}\label{lemma1}
The envelope of the geodesics through $P_0 = (x_0, y_0)$ is the conjugate point locus to $P_0$ for the model metric
$($equation~\eqref{model}$)$ and is the parabola $y = \frac{1}{4 y_0} (x-x_0)^2$ tangent to the boundary $y=0$ at $(x_0, 0)$.
As a~singularity, the envelope realizes the simplest of the stable singularities of maps $\mathbb{R}^2 \to
\mathbb{R}^2$, the fold singularity, whose normal form near $(0,0)$ is $(u,v) \mapsto (u, v^2)$.
$($See~{\rm \cite[\textit{Theorems}~4.4 \textit{and}~4.5]{GG}}  for results on the fold singularity.$)$
\end{lemma}

We prove the lemma in Appendices~\ref{A1} and~\ref{A2}.
In Appendix~\ref{A1} we review the def\/inition of `envelope' and show that the envelope is indeed the conjugate locus.
In Appendix~\ref{A2} we compute our specif\/ic envelope and show that the map for which it is a~singularity is a~simple fold.
More important than the exact formula for the envelope given in the lemma   is the fact that it represents a~stable
singularity.
We also review the def\/inition of the fold and of a~stable singularity in Appendix~\ref{A2}.
We urge the reader to see the discussion in~\cite{Levi}, especially Fig.~5.6
for another good picture and a~discussion of this model example.

\subsection{Higher dimensions}
To place the model example (equation~\eqref{model}) in higher dimensions, take $x = (x_1, \ldots, x_{n-1}) \in
\mathbb{R}^{n-1}$ and work in the upper half space $y\ge 0$ of $\mathbb{R}^n = \mathbb{R}^{n-1} \times \mathbb{R}$ with
coordinates $(x,y)$.
Understand~$dx^2$ to mean the Euclidean metric~$\sum dx_i^2$.
Euclidean rotations about the vertical axes $x = x_0$ are isometries for the model metric.
The conjugate locus is obtained by taking the envelope just worked out in the lemma above for the planar case and
rotating it about the vertical axis through~$P_0$ to obtain a~hypersurface of revolution.

The lemma above holds as is.
In the equation for the conjugate locus we interpret $(x-x_0)^2$ to mean $\|x - x_0\|^2: = \sum (x_i - x_{i 0})^2$.
The singularity is again a~fold.
The normal form for the fold map from $\mathbb{R}^n \to \mathbb{R}^n$ remains the same, remembering to write $(x,y) \in
\mathbb{R}^{n-1} \times \mathbb{R} = \mathbb{R}^{n}$.

\section{Reduction to the model example}

\begin{proof}[Proof of Theorem~\ref{Structure}, the structure theorem]
Scale Seifert's cylinder coordinates by
$(x,y) $ $\to (\epsilon x, \epsilon y)$.
Here $x \in \mathbb{R}^{n-1}$, $y \in \mathbb{R}$, $y \ge 0$ and $\epsilon > 0$ is suf\/f\/iciently small.
We ask the reader to take a~glance at Seifert's metric normal form (equation~\eqref{SeifertModel})
and the function~$f$ there.
Taylor expand~$f$ about the origin: $f(x,y)= 1 + ax + by +O(x^2 + y^2)$ and write $f_1 (x,y) = ax + by$ for the linear
term.
In the rescaled coordinates Seifert's metric normal form (equation~\eqref{SeifertModel}) becomes:
\begin{gather*}
ds^2_{\rm JM} = \epsilon^3 y\Big(dy^2 + \Big(1 + \epsilon f_1 + O\big(\epsilon^2\big)\Big(\sum dx_i^2 + \sum \epsilon^2 h_{ij} (x,y) dx_i dx_j\Big)\Big)\Big)
\\
\phantom{ds^2_{\rm JM}}
= \epsilon^3 y \big(dy^2 + dx^2 + \epsilon f_1 dx^2 + O\big(\epsilon^2\big)\big)
\end{gather*}
where the $O(\epsilon^2)$ term only contains $dx_i$, $dx_j$ terms (no $dy$'s).
Dividing a~metric by a~positive constant does not change its geodesics.
Divide our metric by $\epsilon^3$ to get the metric:
\begin{gather*}
ds^2 = y\big(dy^2 + dx^2 + \epsilon f_1 dx^2 + O\big(\epsilon^2\big)\big),
\end{gather*}
whose conjugate locus is identical (after rescaling) to that of the original Seifert form.

Now view this expression for $ds^2$ as an instance of the Jacobi--Maupertuis principle.
In other words take the energy~$E$ to be $0$ so that the overall conformal factor~$y$ corresponds to the same potential
$V= -({1/2})y$ as in our model example.
View the term in parenthesis as the ``underlying metric''.
Now play the JM game in reverse, to write out Newton's equations, in Hamiltonian form, based on the structure of this
metric.
The kinetic energy metric part of our model has changed from $dy^2 + dx^2$ to $dy^2 + (1 + \epsilon f_1) dx^2 +
O(\epsilon^2)$ where the $O(\epsilon^2)$ error term does not involve $dy$ but only $dx_i dx_j$ terms.
Set $y = x_0$ momentarily so that we can write the metric tensor of this metric in the uniform manner $\sum g_{ab}
dx_a dx_b$ with $a$, $b$ now running from $0$ to $n-1$.
The Hamiltonian whose Hamilton's equations are Newton's equations is $\frac{1}{2} (\sum g^{ab} p_a p_b - y)$ where
$g^{ab}$ is the inverse to the matrix of metric coef\/f\/icients $g_{ab}$.
We see that
\begin{gather*}
(g)_{ab} = \left(
\begin{matrix}
1 & 0
\\
0 & I + \epsilon f_1 I + O\big(\epsilon^2\big)
\end{matrix}
\right)
\end{gather*}
from which it follows that the inverse matrix is
\begin{gather*}
g^{ab} = \left(
\begin{matrix}
1 & 0
\\
0 & I - \epsilon f_1 I + O\big(\epsilon^2\big)
\end{matrix}
\right)
\end{gather*}
yielding the Hamiltonian
\begin{gather*}
H_{\epsilon} = \frac{1}{2} \big(p_x^2 + p_y^2 - \epsilon f_1 p_x^2 + O\big(\epsilon^2\big)\big) - \frac{1}{2}y,
\end{gather*}
where $p_x^2$ means $\sum p_i^2$.
This Hamiltonian is a~small order~$\epsilon$ perturbation of the Hamiltonian
\begin{gather*}
H_0 = \frac{1}{2} \big(p_x^2 + p_y^2\big) - \frac{1}{2}y
\end{gather*}
for our model problem, solved in the last section.

Write $(x,y, p) \mapsto \Phi_t^{\epsilon} (x,y, p)$ for the Hamiltonian f\/low of our perturbed Hamiltonian $H_{\epsilon}$
with the unperturbed model f\/low being $\Phi_t^0$.
Since $H_{\epsilon}$ is within $C \epsilon$ of $H^0$ in the $C^k$ topology over compact sets (any~$k$ up to the
smoothness of the original problem), we have that their associated f\/lows are also close, provided we restrict to compact
subsets.
In other words, if $(x,y, p; t)$ are conf\/ined to vary over a~compact subset of $K \subset \mathbb{R}^n \times
\mathbb{R}^n \times \mathbb{R}$, then the restrictions of the associated Hamiltonian f\/lows $\Phi_t^{\epsilon}$ and
$\Phi_{t}^0$ are $O(\epsilon)$ close in the $C^{k}$ topology.
(Yes,~$k$, not $k-1$.
We get $k \to k -1$ when we dif\/ferentiate~$H$ to get the Hamiltonian vector f\/ield.
But we add $1$ back when we integrate the vector f\/ield to get the f\/low.) Now we worked out the details of the
unperturbed f\/low in in the last section (equation~\eqref{Solutions}).

Why can we restrict to compact sets? The geodesics for the model problem all leave a~given cylinder set in a~bounded
time $T_0$, and hence the same is true of the perturbed problem, with a~perhaps somewhat bigger $T_0$, say $2T_0$ to be
safe.
Fix the point $q = (x_0, y_0)$ (in rescaled variables).
The perturbed geodesics through $q = (x_0, y_0)$ are obtained by solving Hamilton's equations for $p =(p_x, p_y)$ lying
in the sphere $H_{\epsilon} (x_0, y_0, p) = 0$.
Thus, in computing the geodesics and conjugate locus we need only vary~$p$ and~$t$ over a~compact set of the form
$S^{n-1} \times [0, 2T_0]$.

We want to compare the singular loci of the map $(p, t) \mapsto \pi(\Phi_t^{\epsilon} (x_0,y_0, p))$ to that of the
unperturbed map $(p,t) \mapsto \pi(\Phi_t^{0} (x_0,y_0, p))$.
Here $\pi(x,y,p) = (x,y)$ is the projection onto conf\/iguration space.
The unperturbed map $(p,t) \mapsto \pi(\Phi_t^{0} (x_0,y_0, p))$ is the subject of Lemma~\ref{lemma1} of the previous section and
is structurally stable.
Hence there is an $\epsilon_0 >0 $ suf\/f\/iciently small so that for all $\epsilon < \epsilon_0$ the singularities of $(p,
t) \mapsto \pi(\Phi_t^{\epsilon} (x,y, p))$ are all folds (and are close to those of $ \pi \circ \Phi_t^0$).

The f\/lows being within order~$\epsilon$ we know that the maps $(p, t) \mapsto \Phi_t^{\epsilon} (x,y, p)$ and $(p,t)
\mapsto \Phi_t^0 (x,y, p)$ are $C^2$ close for~$\epsilon$ suf\/f\/iciently small.
Take~$\epsilon$ small enough that structural stability holds: the singular locus of the perturbed map is a~fold.
This locus is our conjugate locus.

The conjugate locus must touch the Hill boundary at the brake orbit through~$q$ as before.
Since the only geodesic through~$q$ touching the Hill boundary is the brake orbit, this is the only point where $C(q)$
intersects the boundary.
Since the conjugate locus and the boundary are both smooth hypersurfaces, and $C(q)$ lies entirely on one side of the
boundary, it must touch it tangentially.
\end{proof}

\section{Proof of the Seifert properties (3)--(5) and Theorem~\ref{Theorem3}}\label{Section5}

We begin with the model problem, taking $z = h_0 - gy$ so we can use the ball-throwing analogy.
The steeper the angle of the throw, the closer we get to the Hill boundary $z =h_0$ and only the straight-up throw
touches the boundary.
The cut-of\/f angle of 45 degrees is angle of maximal horizontal throw: at a~f\/ixed speed this is the upward angle to throw
a~ball so as to achieve the maximum horizontal distance before the ball hits the ground again at $z = z_0$.
Any higher angle and the ball drops short of the 45 degree throw, and the corresponding arc hits the conjugate locus
{\it before the ball hits the ground}, i.e.\ closer to the Hill boundary then when we started.
(Any lower angle and the ball hits the ground before it hits the conjugate locus, and hits the ground short of the 45
degree throw.) Thus if the angle of throw with the vertical is less than 45 degrees then the point where the geodesic
hits the conjugate locus $C(q)$ is closer to the Hill boundary than the starting point~$q$.

\begin{proof}[Proof of Seifert properties (3), (4), and (5).] Properties (3), (4) and (5) above regarding the Seifert coordinates are
easily verif\/ied for the model problem.
The reader can work out precise algebraic relations relating angles of steepness to heights.
The conditions involved in the three properties are open conditions in the $C^2$-topology on curves.
The real problem is an order~$\epsilon$ perturbation of the model problem as measured in the $C^1$-topology on the space
of vector f\/ields.
Consequently the geodesics for the real problem lie within an~$\epsilon$-$C^2$ neighborhood of those of the model
problem.
Consequently these properties continue to hold for the real problem, provided we take~$\epsilon$ small enough.
The precise constants involved will need to be relaxed a~bit.
The smaller we take~$\epsilon$, the closer we are to being able to use the same algebraic relations which the reader may
have worked out in the model problem.

{\it Property $(3)$.} Let us see the details of this argument for property~(3).
In the unperturbed model example, all the geodesics are parabolas.
A~bit of algebra shows that they can be written $y-y_m = \frac{1}{y_m} (x-x_m)^2$ where the vertex of the parabola,
which is the minimum value of~$y$ lies at $(x_m, y_m)$.
One then computes that $|dy/dx| \ge 1$ provided $|x-x_m| \ge \frac{1}{2}y_m$ which is to say $y \ge \frac{5}{4} y_m$.
In other words, once we reach a~height of $\lambda y_m$ or greater, $\lambda > 5/4$ along the parabola, the tangent line
to this parabola is less than 45~degrees
from the vertical.

It follows that if we take any constant $\lambda > 5/4$ then there is an~$\epsilon$ suf\/f\/iciently small, such that any
geodesic which enters into the cylinder neighborhood of height $\epsilon_1 < \epsilon$ will be leaving through the roof
of the cylinder of height $\lambda \epsilon$ and with a~tangent direction to the vertical of angle $45$ degrees or less.
The angle of the tangent with the vertical in the model parabolas decreases monotonically with their height, so the same
is true of the perturbed example, for~$\epsilon$ suf\/f\/iciently small.
We can increase~$\lambda$ so as to guarantee that this angle is, say, $42$ degrees, for example.

The constant in property (5) is verif\/ied by rewriting the model problem with a~constant~$g$: $\ddot y = g$ and observing
that~$g$ corresponds to the length of the force, or gradient of~$V$, and then doing a~bit of algebra and scaling.
\end{proof}

\begin{proof}[Proof of Theorem~\ref{Theorem3}]
To prove Theorem~\ref{Theorem3}, recall, as described a~few paragraphs up, that in the model problem the
downward pointed cone ${\rm DC}(q)$ is the cone of vectors making an angle of 45~degrees with respect to the vertical,
regardless of the initial point $q= (x_0, y_0)$.
The real case is a~small perturbation of the model case with the size of the perturbation tending to zero as we tend to
the boundary.
\end{proof}

\section{Proof of main theorem}
\begin{proof}
Consider the regular point $q_0$ of the boundary together with cylindrical neighborhoods centered on $q_0$ for which the
properties of Seifert hold.

Now by property (3) any geodesic which enters into the cylinder of height $\epsilon_1$ must leave through the roof of
a~cylinder of height $\lambda \epsilon_1$ at an angle closer to 42~degrees to the vertical.
Here~$\lambda$ a~f\/ixed constant, somewhat bigger than~$5/4$.
(We could have taken any degree less than~45 in place of~42~degrees.)
By Theorem~\ref{Theorem3}, for $\epsilon_1$ suf\/f\/iciently small,
these geodesic arcs all have conjugate pairs $q_1$, $q_2$ with $q_1$ being at height $y= \lambda \epsilon_1$ and with
$q_2$ being at a~lower height $ y < \lambda \epsilon_1$.
We refer the reader again to Fig.~\ref{parabola}.

We now simply insist that our geodesics enter the cylinder of height $\epsilon_1$ about $q_0$.
Any such geodesic is of the type described in the previous paragraph.
We are guaranteed our conjugate pair along this geodesic.
\end{proof}

\appendix

\section{Envelopes and conjugate locus}\label{A1}

A~$k$-parameter family of immersed curves on a~manifold~$M$ is a~smooth map $\Gamma: X \times \mathbb{R} \to M$ where, for
each $x \in X$, the curve $c_x (\cdot) = \Gamma(x, \cdot)$ is immersed, and where~$X$ is a~smooth~$k$-dimensional
manifold.
The {\it envelope} of the family is its set of critical values.

According to the chain rule, if $\Psi: X \times \mathbb{R} \to X \times \mathbb{R}$ is any dif\/feomorphism, then the set
of criticial values of~$\Gamma$ and $\Gamma \circ \Psi$ are identical.
A~{\it reparameterization} of the $k$-parameter family is a~particular kind of dif\/feomorphism of the form $\Psi(x,t) =(x,
\tau(x,t))$.
Thus each curve $c_x$ has been reparameterized by a~new parameter $s = \tau(x,t)$.
The envelope (critical {\it values}) of the original family~$\Gamma$ and the reparameterized family $\Gamma \circ
\Psi^{-1}$ are equal.

We apply these considerations to the exponential map of Riemannian geometry.
The conjugate locus of a~point~$q$ is the set of critical points of the exponential map $\exp: T_q M \to M$ based at~$q$.
Recall that $\exp(v) = \gamma(1)$ is the time~1
end point of the unique geodesic $\gamma(t)$ through~$q$ having initial
conditions $\gamma(0) = q$, $\dot \gamma (0) = v$.
One proves that $\gamma(t) = \exp(tv)$.
It follows that upon using polar coordinates for $T_q M$, by writing vectors as $v = s \omega$, with $\omega \in S^{n-1}
\subset T_q M$ a~unit length vector, we can think of the exponential map as an ($n-1$)-parameter family of curves:
$\Gamma: S^{n-1} \times \mathbb{R} \to M$; $\Gamma(\omega, s) = \exp(s\omega)$.
(The parameter~$s$ is arclength.) Thus the conjugate locus of~$q$ is the envelope of this ($n-1$)-parameter family of
geodesics through~$q$.
The discussion on reparameterization invariance above holds, showing that whether we compute the conjugate locus
relative to arclength parameterization~$s$, or by Newtonian-time as af\/forded by Proposition~\ref{JMtoNewton}, we get the
same result for the envelope of curves.
The envelope is the conjugate locus.

{\bf Remark about the origin.} Some thought needs to be applied to the case $s =0$ where polar coordinates break down,
i.e fails to be a~dif\/feomorphism.
But the exponential map is known to be a~dif\/feomorphism at $v =0$, so we exclude the point $q = \exp(0v)$ as being in the
conjugate locus.

\section{Computing the envelope}\label{A2}

\begin{proof}
[Proof of Lemma 1] We have found it helpful to put back in the constant gravitational acce\-le\-ration~$g$ so that our
equations are $\ddot x =0$, $\ddot y = g$ with general solution
\begin{gather*}
x = x_0 + v_1 t,
\qquad
y = y_0 - v_2 t + \frac{g}{2} t^2
\end{gather*}
and energy
\begin{gather*}
E = \frac{1}{2}\big(v_1^2 + v_2^2\big) - gy_0 = 0
\qquad
\text{or}
\qquad
v_1^2 + v_2^2= 2gy_0.
\end{gather*}
We parameterize so the 1-parameter family of solutions passing through $P_0 = (x_0, y_0) $ at time~$t$ by an
angle~$\theta$ according to:
\begin{gather*}
v_1 = \sqrt{2 g y_0} \sin(\theta),
\qquad
v_2 = \sqrt{2 gy_0} \cos (\theta).
\end{gather*}
Subsitute this expression for $v_1$, $v_2$ into the general solution to obtain the explicit one-parameter family $\Gamma
(\theta, t) = (x(\theta, t), y(\theta, t))$.
Compute $dx$ and $dy$ in terms of $\theta$, $t$, $d \theta$, $dt$ to arrive at
\begin{gather*}
dx \wedge dy = t(2g y_0 - v_2 g t) dt \wedge d \theta.
\end{gather*}
It follows that the critical points of the map~$\Gamma$ are def\/ined by $t=0$ and $2g y_0 -v_2 gt = 0$.
We ignore the singularity at $t=0$ as a~coordinate singularity.
See the f\/inal remark of the last appendix.
The critical point locus is $v_2 t = 2y_0$.
(Remember $v_2 = v_2 (\theta)$ as above.) Plugging this relation into the general solution and using $v_1^2 = 2g y_0 -
v_1^2$ we f\/ind, after some algebra, that the envelope of the family, being the~$\Gamma$-image of the set of critical
points, is
\begin{gather*}
y = \frac{1}{4 y_0} (x-x_0)^2,
\end{gather*}
as claimed.

The Jacobian of~$\Gamma$ is
\begin{gather*}
d \Gamma = \left(
\begin{matrix}
v_1 & v_2 t
\\
-v_2 + gt & v_1 t
\end{matrix}
\right).
\end{gather*}
From this expression we verify, after a~bit of algebra, that the tangent to the envelope at a~point $\Gamma(\theta, t)$
is indeed $\dd{\Gamma}{t}$, the tangent to the corresponding curve $t \mapsto \Gamma(\theta, t)$.
The direction of this tangency is $(1, \frac{1}{2 y_0} v_1 t)$ in the $(x,y)$-plane.
Back in the domain of~$\Gamma$ it is the vector $(1,0)$ corresponding to the tangent vector $\dd{}{t}$.
This computation fails at the brake point since the tangent to the brake curve $t \mapsto \Gamma(0,t)$ is zero at the
brake instant.
Special considerations are required there and are supplied in a~few paragraphs.

{\bf Fold.}
The condition that a~smooth map $F: \mathbb{R}^n \to \mathbb{R}^n$ represents a~fold singularity is that the
set of critical points forms a~smooth hypersurface, that along this hypersurface the rank of~$dF$ is~$n-1$, and that the
kernel of~$dF$ is transverse to the tangent space to the hypersurface.
We have already computed the tangent space to the envelope everywhere except at the brake point.
At a~point of the envelope the rows of $d\Gamma$ must be linearly dependent.
(Both rows are nonzero for $t \ne 0$.) We read of\/f from our expression for $d \Gamma$ that the kernel of $d\Gamma$ at
a~point of the envelope is then the span of $(v_2 t, -v_1) = (2y_0, -v_1)$ which is the same as the span of $(1, -v_1/
2y_0)$ This establishes that the map is a~fold away from the brake point $(x_0, 0)$.

To establish that the map~$\Gamma$ is a~fold at the brake point we can Taylor expand its components.
For simplicity, set $g =1/2$ now and $t_b = 2 \sqrt{y_0}$ so that the brake point occurs when $(\theta, t) = (0, t_b)$.
Set $t =t_b + h$ and think of $\theta$, $h$ as small.
We compute
\begin{gather*}
x = 2y_0 \theta + O(\theta h),
\qquad
y = \frac{1}{4} h^2 + y_0 \theta^2 + O\big(\theta^2 h\big)
\end{gather*}
almost the normal form for a~fold.
The kernel of $d \Gamma$ is the $\dd{}{h}$ direction, which is also the $\dd{}{t}$ direction.
The envelope in the $\theta h$ coordinates is given by $h = 0$ plus higher order terms, i.e.~the envelope's tangent at
the brake point is spanned by $\dd{}{\theta}$.
The kernel of $d \Gamma$ is transverse to the tangent to the envelope even at the brake point.
The map is still a~fold at the brake point.
\end{proof}

{\bf On stable maps.} A~singularity for a~map is a~map germ near a~critical point.
The singularity is called ``stable'' if whenever any map $F:M \to N$ realizes this germ (is dif\/feomorphic to it), then
there is an~$\epsilon$-neighborhood about~$F$ in the $C^k$-topology on the space of maps (some $k > 0$) such that any
map~$G$ in this neighborhood is dif\/feomorphic to~$F$ in some neighborhood of the critical point.
That is to say, there are dif\/feomorphisms $\psi: N \to N, \phi: M \to M $ so that $G = \psi \circ F \circ \phi^{-1}$ in
a~neighborhood of the critical point.
The Morse lemma asserts that Morse functions are stable.
The simplest stable singularity of maps from manifolds $M$, $N$ of the same dimension are the fold singularities.
See~\cite{GG} for more on the basics of singularity theory.

\subsection*{Acknowledgements}
I thank Mark Levi and Mikhail Zhitomirskii for helpful e-mail conversations.
I acknowledge NSF grant DMS-1305844 for support.

\pdfbookmark[1]{References}{ref}

\LastPageEnding

\end{document}